




\magnification\magstep1
\baselineskip=18pt
\overfullrule = 0pt
\def\vp{\varepsilon}
\def\n{\noindent}

\centerline{\bf Uniform embeddings, homeomorphisms and quotient}
\centerline{\bf maps between Banach spaces (A short survey)}\medskip
\centerline{Joram Lindenstrauss\footnote*{Supported in part by the
U.S.-Israel Binational Science Foundation.}}
\centerline{The Hebrew University}
\centerline{Jerusalem, Israel}\bigskip

{\bf 1. Introduction}

A well known theorem of Mazur and Ulam [MU] states that if $T$ is an
isometry from a Banach space $X$ onto a Banach space $Y$ such that $T0=0$
then $T$ is linear. Another well known (and deeper) theorem due to Kadec
[K] states that any two separable infinite dimensional Banach spaces are
mutually homeomorphic. Thus while the linear structure of a Banach space is
completely determined by its structure as a metric space, the structure of
a Banach space as a topological space contains no information on the linear
structure.

In the present paper I consider the situation in between those extremes.
Already if we weaken the isometry assumption by just a little and consider
``almost isometries'' we encounter interesting problems and results (e.g.\
the work of F.~John in the context of the theory of elasticity [J] or the
work around the Hyers-Ulam problem, see [Gev] and [LS]). I shall not
discuss this topic here and move somewhat more away from isometry. The
topics of our discussion here will be the structure of Banach spaces as
uniform spaces and the Lipschitz structure of Banach spaces. It turns out
that the study of these topics leads to a rich interplay  between various
areas:\ topology, geometric measure theory, probability, harmonic analysis,
combinatorics and of course the geometry of Banach spaces.

In a book [BL] which is now being written and which will (hopefully) appear
in 1998 there is a detailed study of many aspects of the structure of
uniformly continuous functions and in particular Lipschitz functions on
Banach spaces (e.g.\ extension of functions, differentiability, uniformly
continuous selections, approximation theorems, fixed points etc.). It also
contains a study of the almost isometric topics mentioned above. Here I
shall survey the main results on three topics concerning these
functions\medskip

\n (i)~~Uniform and Lipschitz embeddings of one Banach space into another.

\n (ii)~~Uniform and Lipschitz classification of Banach spaces and their
balls.

\n (iii)~~Uniform and Lipschitz quotient maps. \medskip

I will just state the main
results, explain them and give references to the papers in  which they were
originally proved. For complete proofs, additional results and further
references I refer to the forthcoming book.

It is worthwhile to mention that the embedding problems treated in (i)
above have discrete analogues which lead to the study of natural problems
on finite metric spaces (usually graphs with their obvious metric). These
problems are of a combinatorial nature and are connected to topics in
computer science. This direction is however not discussed here and again I
refer to [BL] for a detailed treatment of this topic.

The theory of Lipschitz and uniformly continuous functions on Banach spaces
has been developing in a slow but rather steady pace over the last 35 years
and by now much is known in this direction. Nevertheless many basic and
natural questions remain unanswered. In the last section of this paper I
present a sample of open problems (not necessarily the central ones) which
are related to the material discussed in the other sections.

Finally let me mention that the recent introductory text on Banach space
theory [HHZ] contains in its last chapter (chapter 12) an introduction
(with proofs) to the subject matter of the present survey.\bigskip

{\bf 2. Embeddings}

Let us start by examining Lipschitz embeddings. The main question here is
the following:\ Assume that there is a Lipschitz embedding $f$ from a
Banach space $X$ into a Banach space $Y$ (i.e.\ $f$ is a Lipschitz
injection and $f^{-1}$ is also Lipschitz on its domain of definition). Does
this imply that $X$ is actually linearly isomorphic to a subspace of $Y$?

The main tool for handling this problem is  differentiation. Let us
recall the definition of the two main types of differentiation.

A map $f$ defined on open set $G$ in a Banach space $X$ into a Banach space
$Y$ is called {\sl G\^ateaux differentiable\/}	at $x_0\in G$ if for every
$u\in
X$ $$\lim_{t\to \infty} (f(x_0+tu) - f(x_0))/t = D_f(x_0)u\leqno (*)$$
exists and $D_f(x_0)$ (= the differential of $f$ at $x_0$) is a bounded
linear operator from $X$ to $Y$.

The map $f$ is said to be {\sl Fr\'echet differentiable\/} at $x_0$ if the
limit
$(*)$ exists uniformly with respect to $u$ in the unit sphere of $X$, or,
alternatively, if
$$f(x_0+v) = f(x_0) + D_f(x_0)v + o(\|v\|)\quad \hbox{as}\quad \|v\|\to
0.$$

Note that if $f$ is a Lipschitz map and $X$ is finite dimensional the
notions of G\^ateaux and Fr\'echet derivatives coincide. If, on the other
hand, $\dim X=\infty$ then there are many natural examples of Lipschitz
maps which are G\^ateaux differentiable at a point without being Fr\'echet
differentiable there (this is the source of a major difficulty in the
area).

It is trivial that if $f$ is a Lipschitz embedding then at every point
$x_0$
where $f$ is G\^ateaux differentiable $D_f(x_0)$ is a linear isomorphism
(into).

We are thus naturally led to the question of existence of a G\^ateaux
derivative. An important notion in this context is that of the Radon
Nikodym Property (RNP in short). A Banach space $Y$ is said to have RNP if
every Lipschitz functions $f\colon \ [0,1]\to Y$ is differentiable almost
everywhere. There are many equivalent definitions of RNP (involving e.g.\
vector measures or extremal structure of convex sets in $Y$) but the one
given above is certainly the most natural in our context. Much is known
about RNP. Obviously a subspace of a space with RNP has RNP and a Banach
space has RNP if all its separable subspaces have RNP and that RNP is an
isomorphism invariant. A result which goes back to Gelfand [Gel] is that a
separable conjugate space has RNP and therefore all reflexive Banach spaces
have RNP. The typical examples of spaces which fail to have RNP are $c_0$
and $L_1$ (in $L_1$ for example consider the function $f$ from [0,1] to
$L_1(0,1)$ defined by $f(t)=$ the characteristic function of the
interval $[0,t]$).

The main theorem on G\^ateaux differentiability is the following (proved
independently at about the same time in [Ar], [Chr] and [Man]).

\proclaim Theorem 1. Let $f$ be a Lipschitz function from a separable
Banach space $X$ into a space $Y$ with RNP. Then $f$ is G\^ateaux
differentiable almost everywhere.

Since there is no natural measure on $X$ the term a.e.\ in the statement of
the theorem needs explanation. Actually each of the 3 papers mentioned
above uses a different notion of a.e.\ (and the notions are definitely not
equivalent) but all will suit us here. For example Christensen calls a
Borel set $A$ a null set if there is a Radon probability measure $\mu$ on
$X$ so that  $\mu(A+x)=0$ for all $x\in X$. It is easy to see that if $\dim
X<\infty$ then $A$ is null in the above sense iff $A$ is of Lesbegue
measure 0. It is also not hard to verify that a countable union of null
sets is a null set and that a null set has empty interior.

An immediate corollary of the theorem is the following statement.

{\sl Assume $X$ is separable and that there is a Lipschitz embedding of it
into a space $Y$ with RNP. Then there is a linear isomorphism from $X$ into
$Y$.}

What happens if $Y$ fails to have RNP? Let us check first the case $Y=c_0$.
The following result was proved by Aharoni [Ah1].

{\sl Every separable Banach space is Lipschitz equivalent to a subset of
$c_0$.}

Recall that $c_0$ is a ``small'' space, actually a minimal space in the
following sense. Any infinite dimensional subspace of $c_0$ has in turn a
subspace isomorphic to $c_0$. Thus for example $\ell_p$, $L_p$, $1\le
p<\infty$ or $C(0,1)$ all are not isomorphic to a subspace of $c_0$. Hence
any Lipschitz embedding of such a space into $c_0$ is nowhere G\^ateaux
differentiable. It is also interesting to note that Aharoni's result is
equivalent to the statement that every separable metric space Lipschitz
embeds into $c_0$.

We turn to the other typical example of non RNP space, namely $L_1(0,1)$.
Here the
situation is entirely different. It is very likely that every Banach space
$X$ which Lipschitz embeds (even only uniformly embeds) into $L_1(0,1)$ is
already linearly isomorphic to a subspace of $L_1(0,1)$. This is definitely
the case if $X$ is reflexive. This follows from the discussion following
Theorem~2 below. It is interesting to note that in this case we obtain a
result on linearization of Lipschitz embeddings which apparently cannot be
proved by differentiation.

We turn now to the question of existence of a uniform embedding (i.e.\ a
uniform homeomorphism into) of a Banach space $X$ into a Banach space $Y$.
For this question the only case in which a significant result is known is
the case $Y=\ell_2$. The following theorem was proved by Aharoni, Maurey
and Mityagin [AMM].

\proclaim Theorem 2. A Banach space $X$ is uniformly equivalent to a subset
of $\ell_2$ if and only if $X$ is linearly isomorphic to a subspace of
$L_0[0,1]$.

The space $L_0(0,1)$ is the space of all measurable functions on [0,1] with
the topology of convergence in measure. The space $L_0(0,1)$ itself is of
course not a Banach space and the theorem actually holds for general
topological vector spaces $X$. It is unknown if every Banach space $X$
which is isomorphic to a subspace $L_0(0,1)$ is also isomorphic to a
subspace of $L_1(0,1)$. It is known that this is the case if $X$ is
reflexive. In particular the space $L_p(0,1)$ or $\ell_p$ is uniformly
homeomorphic to a subset of $\ell_2$ if and only if $1\le p\le 2$.

If $T$ is any map from a Banach space $X$ into $\ell_2$ then $K(x,y) =
\langle Tx,Ty\rangle$ is a positive definite kernel of $X$. The proof of
Theorem~2 starts with this observation and then the argument shifts to
examining positive definite kernels on $X$ and much of the argument is
probabilistic in nature. Let us also note that $K(x,y) = e^{-\|x-y\|^2}$ is
a positive definite kernel on Hilbert space. There is a uniform embedding
$T$ of $H$ into itself so that $\langle Tx,Ty\rangle   = K(x,y)$. Since
$K(x,x)=1$ this map $T$ gives a uniform embedding of Hilbert space into its
unit sphere. Thus whenever a Banach space (or a metric space) embeds
uniformly into $\ell_2$ it also embeds uniformly into the unit sphere of
$\ell_2$. In this connection it is of interest to note that from the
results quoted in the next section it follows that there are many examples
of Banach spaces $X$ which do not embed uniformly in $\ell_2$ but whose
unit balls $B(X)$ do embed uniformly in $\ell_2$ (this is the case  e.g.\
for $X=L_p$ or $\ell_p$ with $2<p<\infty$).\bigskip

{\bf 3. Uniform and Lipschitz classification of spaces and balls}

As in the previous section we start with Lipschitz mappings, this time with
a bi-Lipschitz map $f$ between a Banach space $X$ and a Banach
space
$Y$. The natural question which arises is the following:\ is $X$ linearly
isomorphic to $Y$? If
$Y$ is separable (and hence also $X$) and has the RNP we can use Theorem~1
and get a G\^ateaux derivative of $f$. Unfortunately the G\^ateaux derivative
may be an isomorphism into (Ives [I] has constructed a Lipschitz
isomorphism from $\ell_2$ onto itself which has a G\^ateaux derivative at 0,
say, and this derivative is an isomorphism from $\ell_2$ onto a
hyperplane).
So one has to use another method to verify that $X$ is isomorphic to $Y$.
One approach which comes to mind is to be careful with the choice of the
point in which one takes the G\^ateaux derivative (we know that it exists
a.e.\
and in Ives' example only one point behaves badly). It is however unknown
how to do such a choice (one problem here is that it is unknown if a
Lipschitz homeomorphism carries null sets to null sets as it does in finite
dimensional spaces). Another possibility is to use Fr\'echet derivatives.
It is trivial that if $f$ is Fr\'echet differentiable at a point $x_0$ then
$D_f(x_0)$ is a linear isomorphism from $X$ onto $Y$. Actually it is enough
that $f$ is $\vp$-Fr\'echet differentiable at some point for small
enough $\vp$.

We say that $f$ is $\vp$-Fr\'echet differentiable at $x_0$ if there
is a bounded linear operator $T$ from $X$ to $Y$ and a $\delta>0$ so that
$$\|f(x_0+u) - f(x_0) - Tu\| \le \vp\|u\|\quad \hbox{for}\quad \|u\| \le
\delta.\leqno (**)$$
It is trivial to check that if $\vp^{-1}$ is larger than the Lipschitz
constant of $f^{-1}$ and $T$ is given by $(**)$ then $T$ is an isomorphism
from $X$ onto $Y$.

The trouble is that there is no general theorem which ensures existence of
Fr\'echet derivatives (or only $\vp$-Fr\'echet derivatives) in this
situation. The only result on existence of Fr\'echet derivatives of
Lipschitz functions is the following deep result of Preiss [P].

\proclaim Theorem 3. Assume $X$ is a Banach space with separable dual. Then
every Lipschitz function $f$ from $X$ to $R$ is Fr\'echet differentiable on
a dense set.

The assumption that $X^*$ is separable is natural here. If $X=\ell_1$ then
$f(x) = \|x\|$ is nowhere Fr\'echet differentiable and a similar function
can be built on any separable $X$ whose dual is not separable. The trouble
with Theorem~3 is that one has in its conclusion a ``dense set'' and not a
``null set'' in any sense which would allow e.g.\ to find a common point of
Fr\'echet differentiability for any given sequence of functions from $X$ to
$R$ (this would have strong implications related to the problem we consider
in this section). At present it is unknown e.g.\ if every
Lipschitz function from $\ell_2$ to the plane has a point of Fr\'echet
differentiability. For finitely many functions there is however a (also
quite complicated) result on $\vp$-Fr\'echet differentiability
proved in [LP].

\proclaim Theorem 4. Assume that $X$ is a separable superreflexive space,
and let $f$ be a Lipschitz function from $X$ to $R^n$. Then for every
$\vp>0$ $f$ has a point of $\vp$-Fr\'echet differentiability.

A space $X$ is called superreflexive if it has an equivalent uniformly
convex norm (in particular it is reflexive).

If $X$ is $\ell_p, 1<p<\infty$, or more generally a space with an
unconditional basis $\{e_i\}^\infty_{i=1}$ the map $f\colon \ X\to X$
defined by
$$f\left(\sum^\infty_{i=1} \lambda_ie_i\right) = \sum^\infty_{i=1}
|\lambda_i|e_i$$
is a Lipschitz map and is nowhere even 1-Fr\'echet differentiable.
The same map can also be considered as a map from $\ell_r$ to $\ell_s$ if
$s>r$. Thus in all these situations there cannot be any existence theorem
for points of $\vp$-Fr\'echet differentiability for general Lipschitz
functions.

In the next section we shall mention some positive results on 
$\vp$-Fr\'echet differentiability of Lipschitz functions between certain pairs of
infinite-dimensional Banach spaces but these results by their very nature
are not of use for the problem we consider in this section.

Heinrich and Mankiewicz [HM] found however a way to deduce linear
isomorphism from Lipschitz equivalence in some rather general situations.
Their argument is based on differentiation but the linear isomorphism they
find is not a differential of the Lipschitz homeomorphism $f$. It is
constructed from differentials of $f$ and $f^{-1}$ in a rather complicated
way. They showed the following.

{\sl Let $f$ be a bi-Lipschitz map from a conjugate Banach space
$X$
onto a Banach space $Y$. Assume that $f$ is G\^ateaux differentiable at a
point $x_0$. Then $D_f(x_0)$ is an isomorphism of $X$ onto a complemented
subspace of $Y$.}

The complementation assertion here is the new fact. In particular it
follows from this and Theorem~1 that

{\sl If $X$ and $Y$ are Lipschitz equivalent separable conjugate spaces
then each is isomorphic to a complemented subspace of the other.}

In other words $X\approx Y\oplus U$ and $Y\approx X\oplus W$ for some $U$
and $W$. Pelczynski showed that under a mild additional hypothesis this
implies that $X$ is actually isomorphic to $Y$ (it was shown recently by
Gowers that without additional assumptions this is false). Anyhow it
follows
from the previous result that for many concrete pairs of spaces (actually
only one needs to be ``concrete'') Lipschitz equivalence implies linear
isomorphism. In particular we have

\proclaim Theorem 5. If $X$ is $L_p(0,1)$ or $\ell_p$ with $1<p<\infty$ and
if $Y$ is Lipschitz equivalent to $X$ then $Y$ is isomorphic to $X$.

In other words the spaces $L_p(0,1)$ and $\ell_p$, $1<p<\infty$, are
determined by their Lipschitz structure. In [AL1] an example was
constructed of two Banach spaces which are Lipschitz equivalent but
nonisomorphic; one space is $c_0(\Gamma)$ for $\Gamma$ uncountable and the
other is a suitable subspace of $\ell_\infty$. Other examples are known by
now but they are of similar nature. In particular there is no known example
of a pair of separable spaces which are Lipschitz equivalent
but not linearly isomorphic.

We pass now to uniform homeomorphism between Banach spaces. A very simple
but useful fact here is that a uniformly continuous map $f$ defined on a
Banach space (in fact on any metrically convex metric space) is a Lipschitz
map for large distances in the sense that for every $\vp>0$ there is a
$C(\vp)$ so that $\|f(x)-f(y)\| \le C(\vp)\|x-y\|$ whenever $\|x-y\|\ge
\vp$. Thus if $f\colon \ X\to Y$ is a uniformly continuous map
$\lim\limits_{n\to \infty}n^{-1}f(nx)$ is a Lipschitz map if the limit
exists. In general there is no reason to assume that this limit exists but
this can be remedied by passing to ultraproducts.

Recall that if $X$ is a Banach space and ${\cal U}$ is a free ultrafilter on	the integers then $X_{\cal U}$ is defined to be the space of all the
bounded sequences $\tilde x = (x_1,x_2,\ldots)$ of elements of $X$ with
$\|\tilde x\|_{\cal U} = \lim\limits_{n\in {\cal U}} \|x_i\|$ (modulo the
sequences of norm 0) with the obvious vector operations. The space $X$
isometrically embeds into $X_{\cal U}$ via the mapping $x\to
(x,x,x,\ldots)$ but in general $X_{\cal U}$ is much larger than $X$ (if
e.g.\ $X$ is separable and infinite dimensional then $X_{\cal U}$ is
nonseparable). Nevertheless the finite dimensional structure of a space is
not lost by passing to an ultraproduct. Any finite dimensional subspace of
$X_{\cal U}$ is almost isometric to a finite dimensional subspace of $X$.
It is worthwhile to introduce here the notion of Banach Mazur distance
which we will use a few times below. Assume that $E$ and $F$ are finite
dimensional Banach spaces of the same dimension. The Banach Mazur distance
$d(E,F)$ between them is defined to be $\inf\{\|T\|\cdot \|T^{-1}\|\}$ where
the infimum is taken over all linear maps $T$ from $E$ to $F$. (Actually
$\log d(E,F)$ is a proper distance functions but we follow the common
practice of discarding the log). In this terminology what we said on
ultraproducts can be expressed as follows. For every finite dimensional
subspace $E$ of $X_{\cal U}$ and every $\vp>0$ there is a finite
dimensional subspace $F$ of $X$ with $d(E,F) \le 1+\vp$.

If we go beyond the finite dimensional subspaces and ask about the
structure
of $X_{\cal U}$ itself we encounter often very difficult questions. There
are however cases where the situation is known. For example one gets easily
from abstract characterizations of $L_p$ spaces that any ultraproduct of an
$L_p(\mu)$ space (e.g.\ $L_p(0,1)$ or $\ell_p$) is an $L_p(v)$ space for a
different (``huge'') measure space.

Coming back to uniform homeomorphisms one gets that if $f$ is a uniform
homeomorphism from a Banach space $X$ onto a Banach space $Y$ then the map
$\tilde f\colon \ X_{\cal U}\to Y_{\cal U}$ defined by
$$\tilde f(x_1,x_2,\ldots, x_n,\ldots) = (f(x_1), f(2x_2)/2,\ldots,
f(nx_n)/n,\ldots)$$
is a Lipschitz homeomorphism ($\tilde f^{-1}$ is obtained from $f^{-1}$ in
the same manner). In other words.

{\sl Uniformly homeomorphic Banach spaces have Lipschitz equivalent
ultrapowers.}

At the level of the ultrapowers one can now use differentiation (i.e.\
G\^ateaux derivatives) to obtain linear maps. That is particularly useful if
we consider finite dimensional spaces since they are essentially not changed
by passing to an ultrapower. Also by using duality it can be shown that in
this context it is possible to bypass the condition of RNP for existence of
derivatives. In this way one can prove the following result due to Ribe
[Ri1].

\proclaim Theorem 6. Let $X$ and $Y$ be uniformly homeomorphic Banach
spaces. Then there is a $C<\infty$ so that for every finite dimensional
subspace $E$ of $X$ there is a subspace $F$ of $Y$ with $d(E,F)<C$ (and of
course vice versa).

In other words the uniform structure of a Banach space determines (up to a
constant) the linear structure of its finite dimensional subspaces. In
particular if $X$ is $\ell_2$ and $Y$ is uniformly
homeomorphic to $X$ then every finite-dimensional subspace $E$ of $Y$
satisfies $d(E,\ell^n_2)\le C$ for a suitable $n$ and a $C$ independent of
$n$. This trivially leads to the following result due to Enflo [E2].

{\sl A Banach space which is uniformly homeomorphic to $\ell_2$ is already
linearly isomorphic to $\ell_2$.}

This approach works only for $\ell_2$. The space $\ell_2$ is probably the
only separable Banach space whose global structure is determined (up to
isomorphism) by the structure of its collection of finite dimensional
subspaces (though this it is not known as yet there are strong partial
results in this direction which show in particular that all other common
spaces fail to have this property).

What is the situation for other Banach spaces? The first question which
comes to mind in this connection is the following:\ Are $\ell_p$ and
$L_p(0,1)$ uniformly homeomorphic for a fixed $p$, $1\le p<\infty$ $p\ne
2$. On the one hand these spaces have common ultraproducts (and in
particular the same finite dimensional structure) but on the other hand are
not isomorphic. It turns out that they are not uniformly homeomorphic. This
was first proved for $p=1$ by Enflo then for $1<p<2$ by Bourgain [Bo] and
finally by Gorelik [Gor] for $2<p<\infty$. The proof of Gorelik is based on
a general principle which in turn is based on Brouwer's fixed point
theorem.

\proclaim Theorem 7 (Gorelik's principle). Let $f$ be a uniform
homeomorphism from
a Banach space $X$ onto a Banach space $Y$. Assume that $f$ carries a ball
centered at the origin and of radius $r$ in a subspace of finite
codimension in $X$ into the $\rho$ neighborhood of a subspace of infinite
codimension in $Y$, then $w(2r)\ge \rho/4$ where $w$ is the
modulus of uniform continuity of $f$.

Using this principle it was proved in [JLS] that several spaces (besides
$\ell_2$) are determined by their uniform structure. In particular.

\proclaim Theorem 8. Any Banach space which is uniformly homeomorphic to
$\ell_p$ $(1<p<\infty)$ is already linearly isomorphic to $\ell_p$.

Probably most common Banach spaces are determined by their uniform
structure but at present the proof of this is not in sight and any class
seems to require a special method. An obvious reason for this difficulty is
that in general Banach spaces are not determined by their uniform
structure. It was proved by Ribe [Ri2] that if $\{p_n\}^\infty_{n=1}$ is a
sequence strictly decreasing to $p>1$ then the spaces $X$ and $X\oplus
\ell_p$ are uniformly homeomorphic but not isomorphic where $X$ is the
direct sum $\left(\sum\limits^\infty_{n=1} \oplus \ell_{p_n}\right)_1$
(i.e.\ the direct sum is taken in the $\ell_1$ norm). In [AL2] the argument
of Ribe was modified so that it works also if 1 is replaced by $s$; $1<s<p$
and
thus one gets even a superreflexive and separable example. A perhaps more
striking example is presented in [JLS]. Let $T_2$ be the ``2 convexified
Tsirelson space''. The precise definition of this space is rather
complicated and not relevant here. What matters is that $T_2$ is a
superreflexive separable space which does not contain a copy of $\ell_2$
but which is ``close to $\ell_2$'' so that any ultrapower of $T_2$ is
isomorphic to $T_2\oplus \ell_2(\Gamma)$ for a suitable uncountable
$\Gamma$. What is proved on $T_2$ (using the method of Ribe as well as a
variant of Theorem~4  the Gorelik principle and other tools) is that any
space uniformly homeomorphic to $T_2$ is linearly isomorphic to either
$T_2$ or $T_2\oplus \ell_2$ and that $T_2$ and $T_2\oplus \ell_2$ are
uniformly homeomorphic but not isomorphic. In other words the spaces
uniformly homeomorphic to $T_2$ represent exactly two isomorphism classes
of Banach spaces.

The construction of a uniform homeomorphism between nonisomorphic spaces is
based in all the examples mentioned above on the fact that unit balls in
completely different spaces can be mutually uniformly homeomorphic. The
uniform structure of a ball in a Banach space contains in it some
information on the linear structure of the space but as we shall see below
only very little information. From the technical point of view the source
of this difference between balls and the entire space is that for studying
uniformly continuous maps on balls we cannot ``go to infinity'' and
transfer the study to Lipschitz maps on ultraproducts.

We pass now to the study of balls. It was already noted by Mazur [Maz] that
the natural nonlinear map $\varphi_{r,s}$ from $L_r(\mu)$ to $L_s(\mu)$
$(1\le r,s <\infty)$ defined by $\varphi_{r,s}(f) = |f|^{r/s} \hbox{ sign }
f$ is a uniform homeomorphism between the unit balls of these spaces. Thus
the unit balls of separable $L_p(\mu)$ spaces $1\le p<\infty$ are uniformly
homeomorphic to the unit ball of $\ell_2$. This is the situation for a much
larger class of Banach lattices. It turns out that the only obstruction to
the uniform homeomorphism of balls in Banach lattices to the unit ball of
Hilbert space is the presence of large cubes. Enflo [E1] proved that it is
impossible to embed the unit balls of $\ell^n_\infty$ into Hilbert space
with a distortion (i.e.\ $\|f\|_{\rm Lip} \|f^{-1}\|_{\rm Lip}$ where $f$
is the embedding) bounded by a constant independent on $n$. Actually, he
proved a somewhat stronger statement and this is one of the first results in	the theory of embedding discrete metric spaces into Banach spaces to which
we hinted at the end of the introduction. The following theorem was proved
by Odell and Schlumprecht [OS] for discrete Banach lattices (i.e.\ spaces
with an unconditional basis) and by Chaatit [Cha] for general lattices.

\proclaim Theorem 9.  The unit ball of a separable infinite dimensional
Banach lattice $X$ is uniformly homeomorphic to the unit ball of $\ell_2$
if and only if $X$ does not contain finite dimensional subspaces
$\{E_n\}^\infty_{n=1}$ with $\sup\limits_n d(E_n,\ell^n_\infty)<\infty$.

It is interesting to recall the context in which Odell and Schlumprecht
proved Theorem~9. They were interested in the following long standing open
problem concerning Lipschitz (or in this context, equivalently, uniformly
continuous) functions from the unit sphere $\{x_j\|x\|=1\}$ of $\ell_2$ to
the real line. Given such a function $f$ and given $\vp>0$, does there
exist an infinite dimensional subspace $Y$ of $\ell_2$ so that the
restriction of $f$ to the unit sphere of $Y$ is constant up to $\vp$, i.e.\
has an oscillation less than $\vp$ (the so-called ``distortion problem'')?
It is known that there are always finite-dimensional subspaces of $\ell_2$
with arbitrarily large dimension which have such a property. In [OS] this
problem was solved in the negative. There is a Lipschitz function $f$ from
the unit sphere of $\ell_2$ to $R$ so that its restriction to the unit
sphere of every infinite dimensional subspace of $\ell_2$ has an
oscillation larger than $\vp_0$ for some $\vp_0>0$. The construction of $f$
is not explicit (and till now no explicit example is known). They first
work on a Tsirelson type space (like the space $T_2$ mentioned above) and
then transfer the result to $\ell_2$ via a uniform homeomorphism of the
unit sphere of $T_2$ with the one in $\ell_2$ (the transfer is not
automatic though, since a uniform homeomorphism does not carry linear
subspaces to linear subspaces).

There are many more spaces whose unit balls are known to be uniformly
homeomorphic to $B(\ell_2)$ besides lattices. N.J.~Kalton noted that one
can apply the Pelczynski decomposition method to this question and was able
e.g.\ to deduce from this that if $Y$ is an infinite dimensional subspace
of a discrete superreflexive lattice $X$
 then
$B(Y)$ is uniformly homeomorphic to $B(\ell_2)$. The class of  spaces
having this property
can be further extended by using complex interpolation of Banach spaces.
For
instance denote by $C_p$  the Schatten spaces of operators on $\ell_2$;
 it follows by interpolation that $B(C_p)$ is uniformly homeomorphic to
$B(\ell_2)$ for $1<p<\infty$ (it is known that the spaces $C_p$, $p\ne 2$,
do not embed linearly into a superreflexive lattice).

Theorem 9 cannot however be extended to all Banach spaces (i.e.\ the
lattice assumption cannot be simply dropped). Raynoud [Ra] proved that
there is a separable Banach space $X$ so that $B(X)$ is not uniformly
homeomorphic to a subset of $\ell_2$ but $X$ does not contain subspaces
$\{E_n\}^\infty_{n=1}$ with uniformly bounded Banach Mazur distances from
$\ell^n_\infty$. For this space $X$, $B(X)$ is also not uniformly
homeomorphic to $B(c_0)$. The proof of Raynoud	is based on the so-called
``stability theory'' of Krivine and Maurey. \bigskip

{\bf 4. Quotient maps}

The notion dual to embeddings (at least in the linear theory) is that of
quotient maps. For the study of this notion we have first to define
properly the concept of nonlinear (Lipschitz or uniform) quotient map. The
most direct concept which comes to mind -- that of a Lipschitz or uniformly
continuous map from one Banach space onto another turns out not to be the
right one in our context. Bates [Ba] proved the following result.

{\sl For every infinite dimensional Banach space $X$ there is a
continuously Fr\'echet differentiable Lipschitz map $f$ onto any separable
Banach space.}

Thus one has to require more of a Lipschitz quotient map in order to get a
concept which is related to the linear structure of the spaces.

A linear quotient map $T$ from a Banach space $X$ onto a Banach space $Y$
is by the open mapping theorem an open map and there exists a constant
$\lambda>0$ so that $TB_X(x,r)\supset B_Y(Tx,\lambda r)$ for every $x\in X$
and
$r>0$ $(B_X(x,r)$ is the ball in $X$ with center $x$ and radius $r$). It
turns out that this property which is automatic for linear quotient maps
has to be built into the definition of nonlinear quotient maps.

{\sl A Lipschitz map $f$ from a Banach space $X$ onto a Banach space $Y$ is
called a Lipschitz quotient map if there is a $\lambda>0$ so that, for all
$x\in X$ and $r>0$, $f(B_X(x,r)) \supset B_Y(f(x),\lambda r)$.}

{\sl A uniformly continuous map $f$ from $X$ to $Y$ is called a uniform
quotient map if there is a function $\varphi(r)$, with $\varphi(r)>0$ for
every $r>0$, so that $f(B_X(x,r))\supset B_Y(f(x), \varphi(r))$ for all
$r>0$.}

The obvious first question to ask with these definitions is whether the
existence of a uniform (resp.\ Lipschitz) quotient map implies the existence
of a linear one.

We start with Lipschitz quotient maps. In the case of Lipschitz embeddings
G\^ateaux derivatives give a good linearization tool (provided the
derivatives exist; for example if we are in the RNP situation). In the
study
of Lipschitz equivalence G\^ateaux derivatives give considerable information,
though not the complete answer. In the case of quotient maps  G\^ateaux
derivatives may give no information (unless one can find a way to use such
derivatives at ``good'' points and not just at an arbitrary point). In fact
in [BJLPS] it is shown that there is for every $1\le p<\infty$ a Lipschitz
quotient map from $\ell_p$ onto itself whose G\^ateaux derivative, at 0 say,
is identically equal to 0.

On the other hand, if a Lipschitz quotient map is $\vp$-Fr\'echet
differentiable at a point and if $\vp$ is small enough then it is easy to
check that the linear operator appearing in the definition of 
$\vp$-Fr\'echet differentiability $(**)$ must be a linear quotient map. In some
cases one can prove that every Lipschitz function from one Banach space to
another is for every $\vp>0$, $\vp$-Fr\'echet differentiable at some
point.

There are pairs of infinite-dimensional Banach spaces $X$ and $Y$ so that
every linear operator from $X$ into $Y$ is compact. Of course in such a
case $Y$ is not a linear quotient space of $X$. It turns out that in many
situations of this type one can prove the existence of points of
$\vp$-Fr\'echet differentiability for Lipschitz maps from $X$ into
$Y$ and therefore $Y$ is also not a Lipschitz quotient space of $X$.

Here are some specific results of this type proved in [JLPS].

{\sl For the following pairs of spaces $X,Y$ every Lipschitz map from $X$
to $Y$ has for every $\vp>0$ points of $\vp$-Fr\'echet differentiability
and consequently there is no Lipschitz quotient map from $X$ onto $Y$

(i)~~$X$ a $C(K)$ space with $K$ compact countable and $Y$ a space having
RNP.

(ii)~~$X$ has a normalized Schauder basis $\{e_i\}^\infty_{i=1}$ so that
$\left\|\sum\limits_i a_ie_i\right\|\le C\left(\sum\limits_i
|a_i|^r\right)^{1\over r}$ for some $C$ and $r$ and all choices of
$\{a_i\}^\infty_{i=1}$ and $Y$ has a uniformly convex norm with modulus of
convexity $\delta(t)$ satisfying $\delta(t) \ge \lambda t^s$ for some
$\lambda>0$ and $s<r$. (For example $X=\ell_r$, $Y=\ell_s$ with
$r>s\ge2$).}

We just mention that there are also some other cases where properties
weaker than (but related to) $\vp$-Fr\'echet differentiability ensure the
nonexistence of a Lipschitz quotient map from $X$ onto $Y$.

We pass now to results on Lipschitz quotient maps which are valid for
uniform quotient maps as well. Therefore we turn to the discussion of
uniform
quotient maps. As in the case of homeomorphisms it is not hard to prove the
following.

{\sl Assume that there is a uniform quotient map from a Banach space $X$
onto $Y$. Then there is a Lipschitz quotient map from an ultrapower
$X_{\cal U}$ of $X$ onto an ultrapower $Y_{\cal U}$ of $Y$.}

As in the case of Ribe's results for homeomorphism (Theorem~6 above) one
can get a quite general linearization theorem for uniform quotient maps
when one goes down to the ``local level'' i.e.\ to the finite dimensional
setting. The following theorem is proved  in [BJLPS].

\proclaim Theorem 10. Assume that $X$ is a superreflexive Banach space and
there is a uniform quotient map from $X$ onto a Banach space $Y$. Then
there is a constant $C<\infty$ so that for every finite dimensional
subspace $E$ of $Y^*$ there is a finite dimensional subspace $F$ of $X^*$
with $d(E,F)<C$.

This theorem could be deduced from Theorem 4 (combined with the
observation above on ultrapowers). However it suffices to use the following
approximation theorem for Lipschitz functions by affine functions whose
proof is considerably simpler than that of Theorem~4. The proof of this
result does not use differentiation and in fact the part of this statement
concerning the estimate is false if we want to use a derivative as an
affine approximant (even if $X$ itself is the real line).

{\sl Let $f$ be a Lipschitz function defined on a ball $B$ in a
superreflexive space $X$ into a finite dimensional space $Y$. Then for
every $\vp>0$ there is
a ball $B_1$ with radius $\rho$ contained in $B$ and an affine function $g$
on $B_1$ so that $|g(x)-f(x)| \le \vp \rho$ for all $x\in B_1$. Moreover
$\rho$ can be estimated from below in terms of $X,\vp$, the radius of $B$,
$\dim Y$ and the Lipschitz constant of $f$.}

An immediate consequence of Theorem 10 and known facts from the linear
theory of Banach spaces is.

\proclaim Theorem 11. Assume that $Y$ is a uniform quotient of $L_p(0,1)$,
$1<p<\infty$. Then $Y$ is isomorphic to a linear quotient of $L_p(0,1)$.
In particular every uniform quotient of a Hilbert space is isomorphic to a
Hilbert space.

For obtaining further results on Lipschitz or uniform quotient maps it
would be useful to have for quotient maps results in the spirit of
Gorelik's principle. Unfortunately at least in the setting of uniform
quotient maps there seem to be no such results. It is shown in [BJLPS] that
for every $1\le p<\infty$ there is a uniform quotient map $f$ from $\ell_p$
onto $\ell_p$ that maps the unit ball of a hyperplane of $\ell_p$ to the
origin. The Gorelik principle (Theorem~7) shows that such a map $f$ cannot
be represented as a composition of a uniform homeomorphism with a linear
quotient map (in any of the two possible orders).\bigskip

{\bf 5. Some open problems}

(1)~~Which separable Banach spaces $X$ are uniformly homeomorphic to
bounded subsets of themselves?

We mentioned in section 2 that this is the case for $X=\ell_2$. Hence by
Theorem~2 and Mazurs map the same is true for $X=L_p(\mu)$ if $1\le p\le
2$. For the same reason this is false for $X=L_p(\mu)$ with $2<p<\infty$.
Aharoni [Ah2] showed that this is true for $X=c_0$ and therefore for any
Banach space containing $c_0$ (like separable $C(K)$ spaces).

(2)~~Can one characterize the Banach spaces $X$ which uniformly embed into
a fixed space $Y$ in situations which are not immediate consequences of
Theorem~2? In particular which Banach spaces embed uniformly into
$L_p(0,1)$ for some fixed $p, 2<p<\infty$?

(3)~~Assume that $X$ and $Y$ are separable Lipschitz equivalent Banach
spaces. Are they linearly isomorphic?

(4)~~Are the spaces $L_p, 1<p<\infty$, $p\ne 2$ determined by their uniform
structure in the sense that any space uniformly homeomorphic to them is
already isomorphic to them? What about $c_0$ or $\ell_1$?

In the case of $L_p$ it is known that any Banach space uniformly
homeomorphic to $L_p$ is linearly isomorphic to a complemented subspace of
$L_p$. In the case of $c_0$ also much is known (the space must be an
${\cal L}_\infty$ space and cannot have $C(w^w)$ as a quotient
space). In the case of $\ell_1$ nothing is known besides Theorem~6 and the
fact that $\ell_1$ and $L_1$ are not uniformly homeomorphic.

(5)~~Are there other Gorelik-like results on uniform homeomorphisms? For
example can a Lipschitz homeomorphism of $\ell_2$ onto itself map the
Hilbert cube $\{x = (\lambda_1,\lambda_2,\ldots)$; $|\lambda_n|\le 1/n\}$
into a hyperplane?

(6)~~Assume that $X$ is a superreflexive separable Banach space. Is $B(X)$
uniformly homeomorphic to $B(\ell_2)$?

(7)~~Assume that $X$ is a Banach space such that $B(X)$ is uniformly
homeomorphic to a subset of $\ell_2$. Is then $B(X)$ uniformly homeomorphic
to $B(\ell_2)$?

(8)~~Is there a Lipschitz quotient map from $\ell_\infty$ onto $c_0$?

It is known (and easy) that there is a retraction from $\ell_\infty$ onto
$c_0$ which is a Lipschitz map. However this map is far from being a
Lipschitz quotient map.

(9)~~Is every Banach space which is a uniform quotient space of $\ell_p$,
$1<p<\infty$, $p\ne 2$ isomorphic to a linear quotient space of $\ell_p$?

By Theorem 11 such a space is a linear quotient space of $L_p$.

(10)~~Is the Gorelik principle true for Lipschitz quotient maps? More
specifically assume that $f$ is a Lipschitz quotient map from an infinite
dimensional Banach space $X$ onto an infinite dimensional Banach space $Y$.
Is it possible that $f$ maps a ball in a finite codimensional subspace of
$X$ to a single point in $Y$?
\vfill\eject

\centerline{\bf References}

{\parindent=45pt
\item{[Ah1]} I.\ Aharoni, Every separable Banach space is Lipschitz
equivalent to a subset of $c_0$, Israel J.\ Math. 19 (1974), 284--291.

\item{[Ah2]} I.\ Aharoni, Lipschitz maps and uniformly continuous maps
between Banach spaces, Thesis, Hebrew University 1978 (in Hebrew).

\item{[AL1]} I.\ Aharoni and J.\ Lindenstrauss, Uniform equivalence between
Banach spaces, Bull.\ Amer.\ Math.\ Soc. 84 (1978), 281--283.

\item{[AL2]} I.\ Aharoni and J.\ Lindenstrauss, An extension of a result of
Ribe, Israel J.\ Math. 52 (1985), 59--63.

\item{[AMM]} I.\ Aharoni, B.\ Maurey and B.S.\ Mityagin, Uniform embedding
of metric spaces and of Banach spaces into Hilbert spaces, Israel J.\ Math.
52 (1985), 251--265.

\item{[Ar]} N.\ Aronszajn, Differentiability of Lipschitzian mappings
between Fr\'echet spaces, Studia Math. 57 (1976), 147--190.

\item{[Ba]} S.M.\ Bates, On smooth nonlinear surjections of Banach spaces,
Israel J.\ Math. (1997).

\item{[BJLPS]} S.M.\ Bates, W.B.\ Johnson, J.\ Lindenstrauss, D.\ Preiss
and G.\ Schechtman, Nonlinear quotient maps, to appear.

\item{[BL]} Y.\ Benyamini and J.\ Lindenstrauss, Geometric nonlinear
functional analysis, A book in preparation.

\item{[Bo]} J.\ Bourgain, Remarks on the extension of Lipschitz maps
defined on discrete sets and uniform homeomorphisms, Springer Lecture Notes
1267 (1987), 157--167.

\item{[Cha]} F.\ Chaatit, On the uniform homeomorphism of the unit spheres
of certain Banach lattices, Pacific J.\ Math. 168 (1995), 11--31.

\item{[Chr]} J.P.R.\ Christensen, Measure theoretic zero sets in infinite
dimensional spaces and applications to differentiability of Lipschitz
mappings, Coll.\ Anal.\ Funct.\ Bordeaux (1973), 29--39.

\item{[E1]} P.\ Enflo, On a problem of Smirnov, Ark.\ Math. 8 (1969),
107--109.

\item{[E2]} P.\ Enflo, Uniform structures and square roots in topological
groups II, Israel J.\ Math. 8 (1970), 253--272.

\item{[Gel]} I.\ Gelfand, Abstrakte Functionen und lineare Operaloren,
Mat.\ Sb. (N.S.) 4(46) (1938), 235--286.

\item{[Gev]} J.\ Gevirtz, Stability of isometries on Banach spaces, Proc.\
Amer.\ Math.\ Soc. 89 (1983), 633--636.

\item{[Gor]} E.\ Gorelik, The uniform nonequivalence of $L_p$ and $\ell_p$
Israel J.\ Math. 87 (1994), 1--8.

\item{[HHZ]} P.\ Habala, P.\ Hajek and V.\ Zizler, Introduction to Banach
spaces (2 volumes) Matfyzpress, Charles University Prague 1996.

\item{[HM]} S.\ Heinrich and P.\ Mankiewicz, Applications of ultrapowers to
the uniform and Lipschitz classification of Banach spaces, Studia Math. 73
(1982), 225--251.

\item{[I]} D.J.\ Ives, Thesis University College London (in preparation).

\item{[J]} F.\ John, Collected papers (2 volumes) Birkhauser 1985.

\item{[JLPS]} W.B.\ Johnson, J.\ Lindenstrauss, D.\ Preiss and G.\
Schechtman, in preparation.

\item{[JLS]} W.B.\ Johnson, J.\ Lindenstrauss and G.\ Schechtman, Banach
spaces determined by their uniform structure, Geometric and Functional
Analysis 6 (1996), 430--470.

\item{[K]} M.I.\ Kadec, A proof of the topological equivalence of all
separable infinite dimensional Banach spaces, Functional Anal.\ App. 1
(1967), 53--62 (translated from Russian).

\item{[LP]} J.\ Lindenstrauss and D.\ Preiss, Almost Fr\'echet
differentiability of finitely many Lipschitz functions, Mathematika 43
(1996).

\item{[LS]} J.\ Lindenstrauss and A.\ Szankowski, Nonlinear perturbations
of isometries, Asterique 131 (1985), 357--371.

\item{[Man]} P.\ Mankiewicz, On differentiability of Lipschitz mappings in
Fr\'echet spaces, Studia Math. 45 (1973), 15--29.

\item{[Maz]} S.\ Mazur, Une remarque sur l'hom\'eomorphism des champs
fonctionnels, Studia Math. 1 (1930), 83--85.

\item{[MU]} S.\ Mazur and S.\ Ulam, Sur les transformations isometrique
d'espaces vectoriels normes, C.R.\ Acad.\ Sci.\ Paris 194 (1932), 946--948.

\item{[OS}] E.\ Odell and Th.\ Schlumprecht, The distortion problem, Acta
Math. 173 (1994), 259--281.

\item{[P]} D.\ Preiss, Differentiability of Lipschitz functions on Banach
space, J.\ Funct.\ Anal. 91 (1990), 312--345.

\item{[Ra]} Y.\ Raynoud, Espaces de Banach superstables, distances stables
et homeomorphismes uniformes, Israel J.\ Math. 44 (1983), 33--52.

\item{[Ri1]} M.\ Ribe, On uniformly homeomorphic normed spaces, Ark.\ Math.
16 (1978), 1--9.

\item{[Ri2]} M.\ Ribe, Existence of separable uniformly homeomorphic
nonisomorphic normed spaces, Israel J.\ Math. 48 (1984), 139--147.\medskip
}

\bye